\magnification=1200


\hsize=14cm    
\vsize=20.5cm       

\parindent=0cm   \parskip=0pt     
\pageno=1 

\def\ind{\hskip 1cm\relax}


\hoffset=15mm    
\voffset=1cm    
 

\ifnum\mag=\magstep1
\hoffset=-0.5cm   
\voffset=-0.5cm   
\fi


\pretolerance=500 \tolerance=1000  \brokenpenalty=5000

\catcode`\@=11

\newcount\secno
\newcount\prmno
\newif\ifnotfound
\newif\iffound

\def\namedef#1{\expandafter\def\csname #1\endcsname}
\def\nameuse#1{\csname #1\endcsname}

\long\def\ifundefined#1#2#3{\expandafter\ifx\csname
  #1\endcsname\relax#2\else#3\fi}
\def\hwrite#1#2{{\let\the=0\edef\next{\write#1{#2}}\next}}

\toksdef\ta=0 \toksdef\tb=2
\long\def\leftappenditem#1\to#2{\ta={\\{#1}}\tb=\expandafter{#2}%
                                \edef#2{\the\ta\the\tb}}
\long\def\rightappenditem#1\to#2{\ta={\\{#1}}\tb=\expandafter{#2}%
                                \edef#2{\the\tb\the\ta}}

\def\lop#1\to#2{\expandafter\lopoff#1\lopoff#1#2}
\long\def\lopoff\\#1#2\lopoff#3#4{\def#4{#1}\def#3{#2}}

\def\ismember#1\of#2{\foundfalse{\let\given=#1%
    \def\\##1{\def\next{##1}%
    \ifx\next\given{\global\foundtrue}\fi}#2}}

\def\section#1{\medbreak
               \global\def\currenvir{section}
               \global\advance\secno by1\global\prmno=0
               {\bf \number\secno. {#1}}}

\def\subsection{\global\def\currenvir{subsection}
                \global\advance\prmno by1
                {\hskip-.5truemm(\number\secno.\number\prmno)}
\hskip 5truemm}

\def\formule{\global\def\currenvir{formule}
                \global\advance\prmno by1
                {\hbox{\rm (\number\secno.\number\prmno)}}}

\def\proclaim#1{\global\advance\prmno by 1
                {\bf #1 \the\secno.\the\prmno.-- }}


\def\ex#1{\medbreak\global\advance\prmno by 1
                {\bf #1 \the\secno.\the\prmno.}}

\def\qq#1{\medbreak\global\advance\prmno by 1
                { #1 \the\prmno)}}

\long\def\th#1 \enonce#2\endth{\global\def\currenvir{th}
   \medbreak\proclaim{#1}{\it #2}\medskip}

\long\def\comment#1\endcomment{}


\def\isinlabellist#1\of#2{\notfoundtrue%
   {\def\given{#1}%
    \def\\##1{\def\next{##1}%
    \lop\next\to\za\lop\next\to\zb%
    \ifx\za\given{\zb\global\notfoundfalse}\fi}#2}%
    \ifnotfound{\immediate\write16%
                 {Warning - [Page \the\pageno] {#1} No reference 
found}}%
                \fi}%
\def\ref#1{\ifx\labellist\empty{\immediate\write16
                 {Warning - No references found at all.}}
               \else{\isinlabellist{#1}\of\labellist}\fi}

\def\newlabel#1#2{\rightappenditem{\\{#1}\\{#2}}\to\labellist}
\def\labellist{}

\def\label#1{%
  \def\given{th}%
  \ifx\given\currenvir%
{\hwrite\lbl{\string\newlabel{#1}{\number\secno.\number\prmno}}}\fi%

\def\given{section}%
  \ifx\given\currenvir%
{\hwrite\lbl{\string\newlabel{#1}{\number\secno}}}\fi%

\def\given{subsection}%
  \ifx\given\currenvir%
{\hwrite\lbl{\string\newlabel{#1}{\number\secno.\number\prmno}}}\fi%

\def\given{formule}%
  \ifx\given\currenvir%
{\hwrite\lbl{\string\newlabel{#1}{\number\secno.\number\prmno}}}\fi%

\def\given{subsubsection}%
  \ifx\given\currenvir%
  {\hwrite\lbl{\string%
   
\newlabel{#1}{\number\secno.\number\subsecno.\number\subsubsecno}}}\fi
  \ignorespaces}

\newwrite\lbl

\def\openall{\openout\lbl=\jobname.lbl}
\def\closeall{\closeout\lbl}

\newread\testfile
\def\lookatfile#1{\openin\testfile=\jobname.#1
    \ifeof\testfile{\immediate\openout\nameuse{#1}\jobname.#1
                    \write\nameuse{#1}{}
                    \immediate\closeout\nameuse{#1}}\fi%
    \immediate\closein\testfile}%

\def\begin{\lookatfile{lbl}
           \input\jobname.lbl
           \openall}
\let\bye\end
\def\end{\closeall\bye}



\font\eightrm=cmr8         \font\eighti=cmmi8
\font\eightsy=cmsy8        \font\eightbf=cmbx8
\font\eighttt=cmtt8        \font\eightit=cmti8
\font\eightsl=cmsl8        \font\sixrm=cmr6
\font\sixi=cmmi6           \font\sixsy=cmsy6
\font\sixbf=cmbx6


\font\tengoth=eufm10       \font\tenbboard=msbm10
\font\eightgoth=eufm8      \font\eightbboard=msbm8
\font\sevengoth=eufm7      \font\sevenbboard=msbm7
\font\sixgoth=eufm6        \font\fivegoth=eufm5

\font\tencyr=wncyr10       
\font\eightcyr=wncyr8      
\font\sevencyr=wncyr7      
\font\sixcyr=wncyr6

\skewchar\eighti='177 \skewchar\sixi='177
\skewchar\eightsy='60 \skewchar\sixsy='60


\newfam\gothfam           \newfam\bboardfam
\newfam\cyrfam

\def\tenpoint{%
  \textfont0=\tenrm \scriptfont0=\sevenrm \scriptscriptfont0=\fiverm
  \def\rm{\fam\z@\tenrm}%
  \textfont1=\teni  \scriptfont1=\seveni  \scriptscriptfont1=\fivei
  \def\oldstyle{\fam\@ne\teni}\let\old=\oldstyle
  \textfont2=\tensy \scriptfont2=\sevensy \scriptscriptfont2=\fivesy
  \textfont\gothfam=\tengoth \scriptfont\gothfam=\sevengoth
  \scriptscriptfont\gothfam=\fivegoth
  \def\goth{\fam\gothfam\tengoth}%
  \textfont\bboardfam=\tenbboard \scriptfont\bboardfam=\sevenbboard
  \scriptscriptfont\bboardfam=\sevenbboard
  \def\bb{\fam\bboardfam\tenbboard}%
 \textfont\cyrfam=\tencyr \scriptfont\cyrfam=\sevencyr
  \scriptscriptfont\cyrfam=\sixcyr
  \def\cyr{\fam\cyrfam\tencyr}%
  \textfont\itfam=\tenit
  \def\it{\fam\itfam\tenit}%
  \textfont\slfam=\tensl
  \def\sl{\fam\slfam\tensl}%
  \textfont\bffam=\tenbf \scriptfont\bffam=\sevenbf
  \scriptscriptfont\bffam=\fivebf
  \def\bf{\fam\bffam\tenbf}%
  \textfont\ttfam=\tentt
  \def\tt{\fam\ttfam\tentt}%
  \abovedisplayskip=12pt plus 3pt minus 9pt
  \belowdisplayskip=\abovedisplayskip
  \abovedisplayshortskip=0pt plus 3pt
  \belowdisplayshortskip=4pt plus 3pt 
  \smallskipamount=3pt plus 1pt minus 1pt
  \medskipamount=6pt plus 2pt minus 2pt
  \bigskipamount=12pt plus 4pt minus 4pt
  \normalbaselineskip=12pt
  \setbox\strutbox=\hbox{\vrule height8.5pt depth3.5pt width0pt}%
  \let\bigf@nt=\tenrm       \let\smallf@nt=\sevenrm
  \normalbaselines\rm}

\def\eightpoint{%
  \textfont0=\eightrm \scriptfont0=\sixrm \scriptscriptfont0=\fiverm
  \def\rm{\fam\z@\eightrm}%
  \textfont1=\eighti  \scriptfont1=\sixi  \scriptscriptfont1=\fivei
  \def\oldstyle{\fam\@ne\eighti}\let\old=\oldstyle
  \textfont2=\eightsy \scriptfont2=\sixsy \scriptscriptfont2=\fivesy
  \textfont\gothfam=\eightgoth \scriptfont\gothfam=\sixgoth
  \scriptscriptfont\gothfam=\fivegoth
  \def\goth{\fam\gothfam\eightgoth}%
  \textfont\cyrfam=\eightcyr \scriptfont\cyrfam=\sixcyr
  \scriptscriptfont\cyrfam=\sixcyr
  \def\cyr{\fam\cyrfam\eightcyr}%
  \textfont\bboardfam=\eightbboard \scriptfont\bboardfam=\sevenbboard
  \scriptscriptfont\bboardfam=\sevenbboard
  \def\bb{\fam\bboardfam}%
  \textfont\itfam=\eightit
  \def\it{\fam\itfam\eightit}%
  \textfont\slfam=\eightsl
  \def\sl{\fam\slfam\eightsl}%
  \textfont\bffam=\eightbf \scriptfont\bffam=\sixbf
  \scriptscriptfont\bffam=\fivebf
  \def\bf{\fam\bffam\eightbf}%
  \textfont\ttfam=\eighttt
  \def\tt{\fam\ttfam\eighttt}%
  \abovedisplayskip=9pt plus 3pt minus 9pt
  \belowdisplayskip=\abovedisplayskip
  \abovedisplayshortskip=0pt plus 3pt
  \belowdisplayshortskip=3pt plus 3pt 
  \smallskipamount=2pt plus 1pt minus 1pt
  \medskipamount=4pt plus 2pt minus 1pt
  \bigskipamount=9pt plus 3pt minus 3pt
  \normalbaselineskip=9pt
  \setbox\strutbox=\hbox{\vrule height7pt depth2pt width0pt}%
  \let\bigf@nt=\eightrm     \let\smallf@nt=\sixrm
  \normalbaselines\rm}

\tenpoint
\def\pc#1{\bigf@nt#1\smallf@nt}         \def\pd#1 {{\pc#1} }


\catcode`\;=\active
\def;{\relax\ifhmode\ifdim\lastskip>\z@\unskip\fi
\kern\fontdimen2  -1.2 \fontdimen3 \string;}

\catcode`\:=\active
\def:{\relax\ifhmode\ifdim\lastskip>\z@\unskip\fi\penalty\@M\ 
\fi\string:}

\catcode`\!=\active
\def!{\relax\ifhmode\ifdim\lastskip>\z@
\unskip\fi\kern\fontdimen2  -1.1 \fontdimen3 
\string!}

\catcode`\?=\active
\def?{\relax\ifhmode\ifdim\lastskip>\z@
\unskip\fi\kern\fontdimen2  -1.1 \fontdimen3 
\string?}

\def\^#1{\if#1i{\accent"5E\i}\else{\accent"5E #1}\fi}
\def\"#1{\if#1i{\accent"7F\i}\else{\accent"7F #1}\fi}

\frenchspacing

\long\def\ctitre#1\endctitre{%
    \ifdim\lastskip<24pt\vskip-\lastskip\bigbreak\bigbreak\fi
  		\vbox{\parindent=0pt\leftskip=0pt plus 1fill
          \rightskip=\leftskip
          \parfillskip=0pt\bf#1\par}
    \bigskip\nobreak}

\let\+=\tabalign

\def\signature#1\endsignature{\vskip 15mm minus 
5mm\rightline{\vtop{#1}}}

\mathcode`A="7041 \mathcode`B="7042 \mathcode`C="7043 \mathcode`D="7044
\mathcode`E="7045 \mathcode`F="7046 \mathcode`G="7047 \mathcode`H="7048
\mathcode`I="7049 \mathcode`J="704A \mathcode`K="704B \mathcode`L="704C
\mathcode`M="704D \mathcode`N="704E \mathcode`O="704F \mathcode`P="7050
\mathcode`Q="7051 \mathcode`R="7052 \mathcode`S="7053 \mathcode`T="7054
\mathcode`U="7055 \mathcode`V="7056 \mathcode`W="7057 \mathcode`X="7058
\mathcode`Y="7059 \mathcode`Z="705A
 
\def\spacedmath#1{\def\packedmath##1${\bgroup\mathsurround=0pt 
##1\egroup$}%
\mathsurround#1 \everymath={\packedmath}\everydisplay={\mathsurround=0pt }}
 
\def\nospacedmath{\mathsurround=0pt \everymath={}\everydisplay={} }


\def\decale#1{\smallbreak\hskip 28pt\llap{#1}\kern 5pt}
\def\decaledecale#1{\smallbreak\hskip 34pt\llap{#1}\kern 5pt}
\def\puce{\smallbreak\hskip 6pt{$\scriptstyle\bullet$}\kern 5pt}

\def\up#1{\raise 1ex\hbox{\smallf@nt#1}}
\def\qed{\raise -2pt\hbox{\vrule\vbox to 10pt{\hrule width 4pt
                 \vfill\hrule}\vrule}}

\def\cqfd{\unskip\penalty 500\quad\vrule height 4pt depth 0pt width 
4pt\medbreak}

\def\build#1_#2^#3{\mathrel{
\mathop{\kern 0pt#1}\limits_{#2}^{#3}}}
\def\date{\the\day\ \ifcase\month\or janvier\or f\'evrier\or mars\or
avril\or mai\or juin\or juillet\or ao\^ut\or septembre\or octobre\or
novembre\or d\'ecembre\fi \ {\old \the\year}}

\def\dateam{\ifcase\month\or January\or February\or March\or
April\or May\or June\or July\or August\or September\or October\or
November\or December\fi \ \the\day ,\ \the\year}

\def\moins{\mathop{\hbox{\vrule height 3pt depth -2pt
width 5pt}}}
\def\crog{{\vrule height 2.57mm depth 0.85mm width 0.3mm}\kern -0.36mm
[}

\def\crod{]\kern -0.4mm{\vrule height 2.57mm depth 0.85mm
width 0.3 mm}}
\def\moins{\mathop{\hbox{\vrule height 3pt depth -2pt
width 5pt}}}

\def\rond{\kern 1pt{\scriptstyle\circ}\kern 1pt}

\def\hfl#1#2{\nospacedmath\smash{\mathop{\hbox to
12mm{\rightarrowfill}}\limits^{\scriptstyle#1}_{\scriptstyle#2}}}

\def\ghfl#1#2{\nospacedmath\smash{\mathop{\hbox to
25mm{\rightarrowfill}}\limits^{\scriptstyle#1}_{\scriptstyle#2}}}

\def\phfl#1#2{\nospacedmath\smash{\mathop{\hbox to
8mm{\rightarrowfill}}\limits^{\scriptstyle#1}_{\scriptstyle#2}}}

\def\pa{\S\kern.15em}
\def\Z{\hbox{{\bf Z}}}

\def\C{\hbox{\bf C}}

\def\Card{\mathop{\rm Card}\nolimits}

\def\cf{{\it cf.\/}}

\def\dra{\ra\kern -3mm\ra}

\def\Ext{\mathop{\rm Ext}\nolimits}

\def\Id{\hbox{\rm Id}}
\def\ie{\hbox{i.e.}}
\def\Im{\mathop{\rm Im}\nolimits}

\def\isom{\simeq}

\def\Ker{\mathop{\rm Ker}\nolimits}
\def\ldra{\lra\kern -3mm\ra}
\def\loc{{\it loc.cit.\/}}
\def\length{\mathop{\rm length}\nolimits}
\def\lra{\longrightarrow}
\def\llra{\nospacedmath\hbox to 10mm{\rightarrowfill}}
\def\lllra{\nospacedmath\hbox to 15mm{\rightarrowfill}}

\def\NS{\mathop{\rm NS}\nolimits}

\def\Pic{\mathop{\rm Pic}\nolimits}

\def\ra{\rightarrow}

\def\Tr{\mathop{\rm Tr}\nolimits}

\def\tx{\kern -1.5pt -}

\def\cc#1{\hfill\kern .7em#1\kern .7em\hfill}

\def\note#1#2{\footnote{\parindent
0.4cm$^#1$}{\vtop{\eightpoint\baselineskip12pt\hsize15.5truecm
\noindent #2}}\parindent 0cm}

\def\og{\leavevmode\raise.3ex
\hbox{$\scriptscriptstyle\langle\!\langle$}}
\def\fg{\leavevmode\raise.3ex
\hbox{$\scriptscriptstyle\,\rangle\!\rangle$}}

\catcode`\@=12

\showboxbreadth=-1  \showboxdepth=-1

\spacedmath{1pt}\parskip=.7mm
\baselineskip=11.1pt
\font\eightrm=cmr10 at 8pt
\overfullrule=0pt
\input amssym.def 
\input amssym
\def\dra{\dashrightarrow} 
\def\tilde{\widetilde}

\def\a{{\alpha}}
\def\b{{\beta}}

\def\o{{\omega}}

\def\cS{{\cal S}}
\def\cA{{\cal A}}
\def\cC{{\cal C}}

\def\cF{{\cal F}}

\def\cI{{\cal I}}
\def\cJ{{\cal J}}

\def\cL{{\cal L}}
\def\cM{{\cal M}}

\def\cO{{\cal O}}

\def\jd{\bar{\cal J}^d{\cal C}}

\def\j2{\bar{\cal J}^2{\cal C}}
\def\cprod_#1^#2{\raise 2pt
\hbox{$\mathrel{\scriptscriptstyle\mathop{\kern0pt\coprod}
\limits_{#1}^{#2}}$}}

\begin
\ctitre {\bf ON THE EULER CHARACTERISTIC OF GENERALIZED KUMMER VARIETIES} 
\endctitre \medskip
\centerline{{\bf Olivier Debarre\note{1}{\baselineskip=3truemm\rm
Partially supported by the European HCM Project ``Al\-ge\-braic Geometry
in Europe'' (AGE), Contract CHRXCT-940557.}}}


\bigskip 

\ind The aim of this note is to apply the Yau-Zaslow-Beauville method
([YZ], [B1]) to compute the Euler characteristic of the generalized Kummer
varieties attached to a complex abelian surface (a calculation also done in
[GS] by different methods). The argument is very geometric: given an ample
line bundle $L$ with $h^0(L)=n$ on an   abelian surface $A$, such that 
each curve  in
$|L|$ is integral, we construct a projective symplectic
$(2n-2)$-dimensional variety
$J^d(A)$ with a Lagrangian fibration
$J^d(A)\ra |L|$ whose fiber over a point corresponding to a smooth curve
$C$ is the kernel of the Albanese map
$J^dC\ra A$. The Yau-Zaslow-Beauville's method shows that the Euler
characteristic of
$J^d(A)$ is $n$ times the number of genus $2$ curves in $|L|$, to wit
$ n^2\sigma(n)$ (where $\sigma(n)=\sum_{m\mid n}m$). The latter  
computation was also done  in [G], where a general conjecture (proved in
[BL] for $K3$ surfaces) expresses these numbers in terms of quasi-modular
forms: if
$N_r^n$ is the number of genus $r+2$ curves in $|L|$ passing through $r$
general points of $A$, one should have
$$\sum_{n\in{\bf N}}N^n_rq^n=\Bigl( \sum_{n\in{\bf N}}
n\sigma(n)q^n\Bigr)^r
\Bigl(\sum_{n\in{\bf N}} n^2\sigma(n) q^n\Bigr) \ .$$ 
\ind Unlike   the case of $K3$ surfaces, none of these varieties $J^d(A)$
seem to be birationally isomorphic to the generalized Kummer  $K_
{n-1}(A)$ introduced by Beauville in [B2] (a symplectic desingularization
of a fiber of the sum morphism $A^{(n)}\ra A$). However, we check that
this is the case when $A$ is a product of elliptic curves and
$d=n-2$. Using a result of Huybrechts, we conclude that $K_{n-1}(A)$ and
$J^{n-2}(A)$ are diffeomorphic hence have the same Euler characteristic $
n^3\sigma(n)$. I would like to thank D. Huybrechts very much for his help
with theorem \ref{huy}.

\section{The symplectic variety $J^d(A)$}

\ind Let $A$ be a complex abelian surface with a polarization $\ell$   of
type $(1,n)$. Assume that each curve with class $\ell$ is integral (this
holds for generic $(A,\ell)$). Let $\hat A$ be the dual abelian surface.
Let
$\phi_\ell:A\ra\hat A$ be the morphism associated with the polarization
$\ell$; there exists a factorization
 $n\Id_A :A\buildrel{\phi_\ell}\over{\lra}\hat A\buildrel{\phi_{\hat
\ell}}\over{\lra}A$, where $\hat\ell$ is a polarization on $\hat A$ of
type $(1,n)$.

\ind We denote by $\Pic^\ell(A)$  the component of the Picard group of $A$
corresponding to line bundles with class
$\ell$, by $\{\ell\}$   the component of the Hilbert scheme that
parametrizes curves in $A$ with class $\ell$, by $\cC\ra \{\ell\}$  the
universal family, and by $\bar{\cJ}\cC\ra \{\ell\}$ the compactified
Picard scheme of this family ([AK]).

\ind The variety $\bar{\cJ}\cC$ splits as a disjoint union
$\cprod_{d\in{\bf Z}}^{}\jd$, where $\jd$ is a projective variety of
dimension $2n+2$, which parameterizes pairs
$(C,\cL)$ where $C$ is a curve on $A$ with class $\ell$ and
 $\cL$ is a torsion free, rank $1$ coherent sheaf   on $C$ of degree $d$
(i.e.\ with
$\chi(\cL)=d+1-g(C)=d-n$). According to Mukai ([M1], ex. 0.5), $\jd$ can
be viewed as a connected component of the moduli space of simple sheaves
$\cL$ on $A$, and therefore is smooth, and admits a (holomorphic)
symplectic structure. There is a natural morphism  
$$\matrix{\a :\jd& \lra&\{\ell\}&\lra&\Pic^\ell(A)\cr 
\ \ \ (C,\cL)&\longmapsto &C &\longmapsto& [\cO_A(C)]\cr}$$ also defined
by $\a(\cL )=\det \cL$.   For each smooth curve $C$ in $\{\ell\}$, the
inclusion
$C\i A$ induces an Abel-Jacobi map $J^dC\ra A$; this defines a rational map
$$\b:\jd\dra A\times \{\ell\}\lra A$$  which is regular since $A$ is an
abelian variety and $\jd$ is normal.  Let
$J^d(A)$ be a fiber of the map $(\a,\b):\jd\lra \Pic^\ell(A)\times A$ (they
are all isomorphic). Note that $J^d(A)$,
$J^{d+2n}(A)$ and
$J^{-d}(A)$ are isomorphic.

\th Proposition
\enonce The symplectic structure on $\jd$ induces a symplectic structure on
the $(2n-2)$\tx dimensional variety $J^d(A)$.
\endth

{\it Proof.}  Recall that there is a canonical isomorphism $T_\cL\jd\isom
\Ext^1(\cL,\cL)$, and that the symplectic form $\o$ is the pairing
$$\Ext^1(\cL,\cL)\otimes\Ext^1(\cL,\cL)\ra
\Ext^2(\cL,\cL)\buildrel{\Tr}\over{\lra}H^2(A,\cO_A)\isom\C
$$
\ind The map
$T_\cL\a$ is the trace map
$T:\Ext^1(\cL,\cL)\ra H^1(A,\cO_A)$, whereas the tangent map at the origin
to the map
$\iota:\Pic^0(A)\ra \jd$ defined by $\iota(P)=P\otimes\cL$ is the dual
$T^*: H^1(A,\cO_A)\ra \Ext^1(\cL,\cL)$. Since $\a\iota$ is constant,
$T\circ T^*=0$; in particular
$$\Ker T \supset \Im T^*=(\Ker T )^\bot\ .$$
\ind Note  also that $\beta\iota\phi_\ell=n\Id_A$ (use the
Morikawa-Matsusaka endomorphism), hence $T_\cL\b\circ
T^*=T\phi_{\hat\ell}$ and $\Ker T_\cL\b\cap
\Im T^*=\{ 0\}$. Since both
$\Ker T$ and $\Ker T_\cL\b$ have codimension $2$, this implies
$$\Ker T=(\Ker T\cap \Ker T_\cL\b)\oplus(\Ker T )^\bot\ ,$$  and the
restriction of $\omega$ to $\Ker T\cap \Ker T_\cL\b=T_\cL J^d(A)$ is
non-degenerate.\cqfd 

\ind The map $\a$ restricts to a   morphism $\a:J^d(A)\ra |L|$ whose fiber
$K^d(C)$ over the point corresponding to a smooth curve $C$ is the
(connected) kernel of the Abel-Jacobi map
$\b:  J^dC\ra A$; it is a Lagrangian fibration.

\section{The Euler characteristic of $J^d(A)$}

\ind We calculate the Euler characteristic of $J^d(A)$ by using the
Lagrangian fibration $\a:J^d(A)\ra |L|$, as in [B1]. 

\th Proposition
\enonce Let $C$ be an integral element of
$|L|$. The Euler characteristic of $K^d(C)$ is
$n$ if the normalization of $C$ has genus
$2$, and $0$ otherwise.
\endth

{\it Proof.}   Let $\eta:\tilde C\ra C$ be the normalization. There is a
commutative diagram (as in
\S 2 of [B1], we may restrict ourselves to the case $d=0$)
  $$\matrix{K(C)&\lra&\bar JC&\buildrel{\b}\over{\lra}&A\cr &&\cup&&\ \
\uparrow \pi\cr &&JC&\buildrel{\eta^*}\over{\lra}&J\tilde C&\ra&0\cr
&&\cup&&\uparrow\cr &&\hskip -6mm(\eta^*)^{-1}(\Ker \pi)&\lra&\Ker
\pi&\ra&0\cr
\cr}$$
\ind By lemma 2.1 of \loc , the group $JC$ acts freely on $\bar JC$. Note
also that for $M$ in $JC$ and $\cL$ in $\bar JC$,
$$\b(M\otimes\cL)=\b(M)+\b(\cL)$$ because this is true when $\cL$ is
invertible, and $JC$ is dense in
$\bar JC$. It follows that
$(\eta^*)^{-1}(\Ker \pi) $ acts (freely) on $K(C)$. As in prop. 2.2  of
\loc , it follows that
$e(K(C))=0$ if $\Ker \pi$ is infinite, that is if $g(\tilde C)>2$.

\ind Assume now that $\tilde C$ has genus $2$. The situation  here is much
simpler than in \loc , because the normalization $\eta$ of $C$ is {\it
unramified}: it is the restriction to $\tilde C$ of the isogeny
$\pi:J\tilde C\ra A$. If $\check C\ra C$ is the minimal unibranch partial
normalization (\cf\ \loc ), it follows   that
$\tilde C\ra\check C$  is an  unramified homeomorphism, hence an
isomorphism (EGA IV, 18.12.6).

\ind There is a commutative   diagram
$$\matrix{0&\lra&\Ker \pi&\lra&J\tilde C&\buildrel{\pi}\over{\lra}&A \cr 
&&\cap&&\ \ \
\cap\,\eta_*&&\vert\vert\cr  0&\lra&K(C) &\lra&\bar J
C&\buildrel{\b}\over{\lra}&A\cr }$$  and an exact sequence
$$1\ra \cO_{\tilde C}^*/\cO_C^*\lra JC\lra J\tilde C\ra 0 \ .$$
\ind If one chooses a line bundle $M$ on $C$ corresponding to a point of
$\cO_{\tilde C}^*/\cO_C^*$ as in the proof of prop. 3.3 of
\loc, it acts on $\bar JC$, hence on 
$K(C)$. Beauville's reasoning proves that $M$ acts {\it freely} on the
complement of 
$ \eta_* J\tilde C$ in $\bar JC$, hence also on the complement of 
$\Ker \pi$ in $K(C)$. It follows that
 $e(K(C))=e(\Ker \pi  )=n$.\cqfd

\ind As a corollary, we get, assuming   that   each curve with class
$\ell$ is integral,
$$e(J^d(A))=n\Card\{\ C\in |L|\ \big| g(\tilde C)=2\ \} \ .$$
\ind 
 It remains to count the number of (integral) genus $2$ curves $C$ in
$|L|$.
 The normalization $\eta:\tilde C \ra C$   induces an isogeny
$\pi:J\tilde C\ra A$ such that $\pi_*\tilde C\in |L|$. Let $r$ be the
degree of
$\pi$; then $\pi^*L$ is numerically equivalent to $r\tilde C$, hence
$rL^2=r^2\tilde C^2=2r^2$ and $r=n$, and $\pi^*\ell$ has type
$(n,n)$. 
  
\ind The number of isomorphism classes of isogenies $\pi:\tilde A\ra A$
such that $\pi^*\ell$  is of type
$(n,n)$  is also the number of isomorphism classes of isogenies
$\hat\pi:\hat A \ra \tilde A$ where
$\tilde A$ has a principal polarization $\theta$ such that
$\hat\pi^*\theta=\hat\ell$, hence also the number of subgroups of
$\Ker\phi_{\hat\ell}$ that are maximal totally isotropic for the Weil form.
This kernel is isomorphic to $(\Z/n\Z)\times(\Z/n\Z)^*$, and the Weil form
is given by
 $e ( (x,x^*),(y,y^*) )=y^*(x)-x^*(y) 
$. Given a quotient group $H$ of $\Z/n\Z$ and any homomorphism $u:H^*\ra 
H$, the set of pairs $(x,x^*)$ in $({\bf Z}/n{\bf Z})\times H^* $ such
that the class of
$x$ in $  H$ is
$u(x^*)$ is such a subgroup, and they are all of this form. Their number
is $$\sum_{ {\bf Z}/n{\bf Z}\twoheadrightarrow H}|H|=\sum_{m\mid
n}m=\sigma(n)\ .$$
\ind To each   isogeny $\pi:J\tilde C\ra A$ correspond  $n^2$ curves in of
genus $2$ in
$|L|$, to wit the curves
$\pi(\tilde C)+x$, for each $x\in \Ker\phi_\ell$. So we get a total of
$n^2\sigma(n)$ such curves, and they are all distinct. 
  
\th\label{euler} Proposition
\enonce Assume each curve with class $\ell$ is integral; then
$e(J^d(A))=n^3\sigma(n)$.
\endth

\section{A degeneration of $J^{n-2}(A)$}

\ind Our aim is to relate the symplectic variety $J^d(A)$ constructed above
with the generalized  Kummer variety $K_{n-1}(A)$. Contrary to the case of
K3 surfaces, these  varieties do not seem to be birational for general
$A$ (except when $n=2$). However, we will prove  that it is the case when
$A$ is a product of elliptic curves and $d=n-2$. For this, we will use the
Mukai-Fourier transform for sheaves on $A$.

\ind For any sheaf $F$ on $A$, we denote by $\cF^\bullet F$ the cohomology
sheaves of the Mukai-Fourier transform of $F$ (see ([M2]). If only $\cF^j
F$ is non-zero, we say that $F$ has weak index $j$, and we write
$\hat F=\cF^jF$; in that case, $\hat F$ has weak index $2-j$, and 
$\hat\cF\hat F\isom(-1)^*F$ (\loc , cor. 2.4). If $H^i(A,F\otimes P_{\hat
x})=0$ for all $\hat x\in\hat A$ and all
$i\ne j$, we say that $F$ has   index $j$; it implies that   $F$ has weak
index
$j$. 

\ind For any $\hat x$ in $\hat A$, we denote by $P_{\hat x}$ the
corresponding line bundle on $A$; we identify the dual of $\hat A$ with
$A$, so that, for any $x$ in $A$, $P_x$ is a line bundle on $\hat A$. 

\ind Let $\cL$ be a sheaf on $A$ corresponding to a point of $\jd$  with
smooth support $C$. For $\cL$ generic in $J^dC$, the surface
$ \cL \otimes \Pic^0(A)$ does not meet the subvariety
$W_d(C)$ of $J^dC$, as soon as $g(C)>2+d$, \ie\ $d<n-1$. In that case, one
has
$H^0(A,\cL\otimes P_{\hat x})=0$ for all $\hat x$ in $\hat A$, so that  
$\cL$ has index $1$, and $\hat\cL$ is a locally free simple sheaf on
$\hat A$ of rank $n-d$, first Chern class $\hat\ell$  and Euler
characteristic $0$.

\th Proposition
\enonce Assume that the N\'eron-Severi group of $A$ is generated by
$\ell$. For $d<n-1$ and
$\cL$ generic in $\jd$, the vector bundle
$\hat\cL$ on $\hat A$ is $\hat\ell$-stable.
\endth

{\it Proof.} We   follow   [FL]: assume $\hat\cL$ is not stable, and look
at torsion-free non-zero quotients of $\hat{\cal L} $ of smallest degree,
and among those, pick one,
$Q$, of
  smallest rank. Because $\NS (\hat A)=\Z\hat\ell$, the degree of $Q$ is
non-positive. The proofs of lemmes  2 and 3 of [FL] apply without
change:   $Q$   has   index
$1$ and  if $K$ be the kernel of $\hat{\cal L}\to Q$, the sheaf $\cF^2 K$
has finite support. Consider the exact sequence 
$$0\ra \cF^1K\ra (-1)^*\cL\ra \hat Q\ra\cF^2K\ra 0\ ;
$$ since $c_1(\hat Q)\cdot\ell=c_1(Q)\cdot\hat\ell\le 0$ ([FL], lemme 1), 
the torsion sheaf
$\hat Q$ has finite support, hence index $0$. But this index is also  $
2-{\rm ind}\, Q=1$; this contradiction proves the proposition.\cqfd

\ind For each $d<n-1$, we have constructed a birational rational map
between
$\jd$ and an irreducible component $\cM^0_{\hat A} (n-d,
\hat\ell,0)$ of the moduli space $\cM_{\hat A} (n-d,
\hat\ell,0)$ of
$\hat\ell$-semi-stable sheaves on
$\hat A$ of rank $n-d$, first Chern class  $ \hat\ell$ and Euler
characteristic
$0$. This map is a morphism if $d<0$.  Let us interpret the maps $\a$ and
$\b$ in this context. Let $(C,\cL)$ be a pair  corresponding to a point of
$\jd$; it follows from the exact  sequence
$0\ra
\cO_C\ra \cO_C(x)\ra \C_x\ra 0$ that
$\det\widehat{\cO_C(x)}\isom\det\widehat{\cO_C}\otimes P_{-x}$, hence
$$\det\hat\cL\isom\det\widehat{\cO_C}\otimes
P_{-\b(\cL)}\isom\det\widehat{\cO_A(-C)}\otimes P_{-\b(\cL)}\ .
$$
\ind Hence, the fibers of $(\a,\b)$ are also the fibers of the map 
$\jd\ra\Pic^\ell(A)\times \Pic^{-\hat\ell}(\hat A)$ which sends $\cL$ to 
$(\det\cL,\det\cF^\bullet\cL)$. Let $ \gamma:\cM^0_{\hat A}  (n-d,
\hat\ell,0)\ra\Pic^{-\ell}(A)\times\Pic^{\hat\ell}(\hat A)$ be the map 
$E\mapsto (\det\cF^\bullet E,\det E)$, and let $M_{n-d}(\hat A)$ be a
fiber. We have proved the following.

\th\label{bira} Proposition
\enonce Assume that the N\'eron-Severi group of $A$ is generated by
$\ell$. For
$d<n-1$, the Fourier-Mukai transform induces a birational isomorphism
between
$\jd$  and an irreducible component
 of $ \cM_{\hat A}  (n-d,
\hat\ell,0)$ which sends $J^d(A)$ onto  $M_{n-d}(\hat A)$.
\endth

\medskip 

\ind We will now study the case where $n-d=2$ and $ A$ is the product of
two general elliptic curves $ F$ and $  G$, with $  \ell$ of bidegree
$(1,n)$. One has  $\hat A=\hat F\times\hat  G$, and
$\hat \ell$ has bidegree
$(n,1)$.  To avoid non-stable semi-stable sheaves, we will study the
moduli space
$\cM'_{\hat A}$ of rank $2$ sheaves
  on
$\hat A$ with first Chern class  $ \hat\ell$ and Euler characteristic
$0$ which are semi-stable for the polarization $\hat\ell'$ of bidegree
$(n+1,1) $, and call $M'(\hat A)$ a fiber of the map
$\gamma:\cM'_{\hat A}\ra\Pic^{-\ell}(A)\times\Pic^{\hat\ell}(\hat A)$
defined above. 

\th\label{fg} Proposition
\enonce The moduli space $\cM'_{\hat A}$ is smooth    and birational to $
\hat A^{(n)}\times A$.   The variety $M'(\hat A)$
  is smooth and birational to  
$K_{n-1}(\hat  A)$.
\endth 

{\it Proof.} Let $E$ be an 
$\hat\ell'$-semi-stable rank $2$ torsion free sheaf on $\hat A$ with first
Chern class  $\hat\ell$ and Euler characteristic $0$. Let $x\in A$; by
semi-stability of
$E^*$, one has
$H^2(\hat A,E\otimes P^{-1}_x)=0$, hence $h^0(\hat A,E\otimes P^{-1}_x)=
h^1(\hat A,E\otimes P^{-1}_x)$. Since $\hat\cF E$ is non-zero, for at
least   one
$x$,
 these numbers are non-zero and there is an inclusion
$P_x\hookrightarrow E$; let $K$ be the kernel of
 $E\ra E/P_x\ra (E/P_x)/(E/P_x)_{\rm tors}$. There is  an exact sequence
$$0\ra K\ra E\ra \cI_Z\otimes K'\ra 0\ ,\leqno{(*)}$$ where   $K'$ is a
line bundle. The line bundle $K$  has bidegree
$(a,b)$, with $a$ and $b$ non-negative and $b(n+1)+a\le (2n+1)/2$ (by
$\hat\ell'$-semi-stability); hence $b=0$ and  
$Z$ is a subscheme of 
$\hat A$ of length
 $n-a $.   

\ind  Set $M=K'\otimes  K^{-1}$. By Serre duality,   $\Ext^1_{\hat
A}(\cI_Z\otimes K',K)$   and $ H^1(\hat A, \cI_Z\otimes M)^*$ are
iso\-mor\-phic. Assume
$H^0(\hat A, \cI_Z\otimes M)=0$, one has $$h^1(\hat A, \cI_Z\otimes
M)=\length (Z)-\chi (\hat A,M)=a $$ and $a>0$   (otherwise $\cI_Z\otimes
K'$ would be a subsheaf of $E$ with $\hat\ell'$-slope $2n+1$),   and
$E$ depends on at most
$2n+3-a$ parameters ($2$ for
$K$,
$2$ for $K'$, $2(n-a)$ for $Z$ and $a-1$ for the extension). Since each
component
 of $\cM'_{\hat A}$ has dimension $2n+2$, this forces
$a=1$ for $E$   generic. Let
$\cM^0$ be the subset of
$\cM'_{\hat A}$ parametrized in this fashion.

\ind Assume now $H^0(\hat A, \cI_Z\otimes M) \ne 0$; one checks (by
projecting onto
 $|M|$), that the set of pairs $(Z,D)$ with $D\in |M|$ and $Z\i D$, has
dimension
$\le n-2a-1+n-a$. Hence $E$ depends on at most $2n-3a-1-\chi(\hat A,
\cI_Z\otimes M)+4=2n-2a+3$ parameters. For $E$   generic, this forces
$a=0$,
$Z$ reduced and
$h^0(\hat A, \cI_Z\otimes M)=1$. This yields a component of $\cM'_{\hat A}$
which can be parametrized as follows. Let $Z=(\hat f_1,\hat
g_1)+\ldots+(\hat f_n,\hat g_n)$ be   generic in $\hat A^{(n)}$, set 
 $L=\cO_{\hat F}(\hat f_1 +\ldots+ \hat f_n )$, and let $ f\in  F$ and
$\hat g\in \hat G$. The vector space  $\Ext^1_{\hat A}( \cI_Z\otimes
p_{\hat F}^*L
\otimes p_{\hat G}^*\cO_{\hat G}(\hat g) ,\cO_{\hat A})$ has dimension
$1$, hence there is a unique extension
$$0\ra p_{\hat F}^*P_f\ra E \ra \cI_Z\otimes p_{\hat F}^*(L\otimes
P_f)\otimes p_{\hat G}^*\cO_{\hat G}(\hat g)\ra 0\ ,$$ where $E$ is
locally free (it satisfies the Cayley-Bacharach condition; see for example
th. 5.1.1 of [HL])
 and stable (the only thing to check is $H^0(\hat A, E\otimes  p_{\hat
F}^*P_{-f}\otimes p_{\hat G}^*\cO_{\hat G}(-\hat g))=0$, and this is true
because the extension is non-trivial). This yields a rational map
$$\phi:\hat A^{(n)}\times F\times \hat G\dra  \cM'_{\hat A}$$ which  is
birational onto its image: given   a locally free
$E$ as above, one recovers $ f$ and the $(\hat g_i-\hat g)$'s   by noting 
that the set $C_E=\{\ x\in A\mid H^0(\hat A,E\otimes P_x)\ne 0\ \}$ is
$$(\{\,- f \,\}\times  G)\cup\bigcup_{i=1}^n( F\times\{\ [\cO_{\hat
G}(\hat g_i-\hat g)]\
\})\ ,$$ the $\hat f_i$'s because $\Ext^1_{\hat A}( \cI_Z\otimes p_{\hat
F}^*L \otimes p_{\hat G}^*\cO_{\hat G}(\hat g) ,\cO_{\hat A})$ must be
non-zero, and $\hat g$ by noting that
$\det E\isom p_{\hat F}^*(L\otimes P_{2f})\otimes p_{\hat G}^*\cO_{\hat
G}(\hat g)$. Because
$H^0(\hat A,E\otimes p_{\hat F}^*\bigl(P_{-f}\otimes\cO_{\hat
F}(-f_1)\bigr)\otimes p_{\hat G}^*\cO_{\hat G}(\hat  g_1-\hat g))$ is
non-zero,
 there exists an exact sequence $(*)$ with $K$ of bidegree
$(1,0)$. This proves that the  set $\cM^0$ defined above is contained in 
the image of
$\phi$, which must therefore be
$\cM'_{\hat A}$.

\ind Finally, $E$ has weak index $1$, $\cF^1E$ has support on $C_E$, and
fixing
$\det\cF^1E$ amounts to fixing
$[\cO_A(C_E)]$. It follows that taking a fiber of $\gamma$
  amounts to fixing
$f$,
$\sum(\hat  g_i-\hat g)$, $\sum \hat  f_i $ and $\hat g$; hence $M'(\hat
A)$ is birational to $K_{n-1}(\hat A)$.\cqfd


 \medskip

\ind The following proof is due to D. Huybrechts, and uses ideas from   
prop. 2.2 of [GH]. 

 \th\label{huy} Theorem
\enonce Let $(A,\ell)$ be a   polarized abelian surface of type $(1,n)$
whose  N\'eron-Severi group  is generated by
$ \ell$. The symplectic varieties $J^{n-2}(A)$, $M_2(\hat A)$ and
$K_{n-1}(\hat A)$ are  de\-for\-ma\-tion equivalent. In particular, they
are all irreducible symplectic.
\endth
 
{\it Proof.} Let $f:\hat\cA\ra S$ be a family of polarized abelian
surfaces, where $S$ is smooth quasi-projective, with
  a relative polarization $\hat\cL$ of type
$(1,n)$, such that the fiber over a point $0\in S$ is   $\hat F\times
\hat G$ with a polarization of bidegree $(n,1)$; assume also that the
N\'eron-Severi group of a very general fiber of $f$ has rank $1$. Let
$g:\cM\ra S$ be  the (projective) relative moduli space of
$\hat\cL$-semi-stable sheaves of rank $2$ with first Chern class $\hat\ell$
and Euler characteristic
$0$ on the fibers of $f$ (\cf\ [HL], th. 4.3.7, p. 92).

 \th Lemma
\enonce Under the hypothesis of the proposition, any rank $2$ torsion
free  sheaf on 
$\hat A$ with first Chern class
$\hat\ell$   which is either simple or semi-stable   is stable.
\endth
 
{\it Proof.} Assume that a sheaf $E$ with these numerical characters is
not stable. There exists an exact sequence
$$0\ra K\ra E\ra \cI_Z\otimes K'\ra 0\ ,$$ where  $K$ and $K'$ are line
bundles on $\hat A$ with $c_1(K)=k\hat\ell$,
$c_1(K')=(1-k)\hat\ell$ and $k>0$. This proves that $E$ is not semi-stable;
moreover, $K\otimes K'^{-1}$ is ample, hence there exists a non-zero
morphism
$u:  K'\ra K $, which induces an endomorphism $E\ra  \cI_Z\otimes
K'\buildrel{u}\over{\ra}K\ra E$ which is not a homothety, and $E$ is not
simple.\cqfd

\ind By the lemma,  the (closed) locus of non-stable points in $\cM$ does
not project  onto $S$. By replacing $S$ with an open subset, we may assume
that there are no such points. Let now $\cS\ra S$ be the (smooth) relative
moduli space of simple sheaves on the fibers of
$f$ (see [AK]). There are embeddings $\cM \i
\cS $ and $\cM'_{\hat F\times
\hat G}\i\cS_0$
 as closed and open subsets. Let $\cS'=\cS\moins(\cS_0\moins \cM'_{\hat 
F\times
\hat G})$; it  is
 open in $\cS$, hence   smooth over $S$. Let
$ \cM'$ be the closure of $g^{-1}(S\moins\{0\})$ in  $\cS'$; the fibers of
$g':\cM'\ra S$ are projective off $0$, and contained in $\cM'_{\hat F\times
\hat G}$ over $0$. Norton's criterion ([N]) shows that points in
$ \cM'_0$ are separated in the moduli space of simple sheaves on $\hat\cA$
(because they are stable), hence in $\cS$; therefore, $ \cM'$ is separated.
By semi-continuity, $ \cM'_0$ is a closed subset of $\cM'_{\hat F\times
\hat G}$ of the same dimension, hence they are equal. Using the lemma, we
get, after shrinking
$S$ again,   a  proper family $g':\cM'\ra S$ with projective irreducible
smooth fibers which coincide with $g:\cM\ra S$ off $0$.

\ind By prop. \ref{bira}, $J^{n-2}(A)$ is birationally isomorphic to  
$M_2({\hat F\times\hat G})$, and   we just saw that the latter deforms to 
$M'({\hat A})$, itself birationally isomorphic to $K_{n-1}(\hat F\times\hat
G)$ by prop. \ref{fg}; in particular, these symplectic varieties are all
irreducible symplectic. Since  birationally isomorphic smooth projective
irreducible symplectic varieties are deformation equivalent ([H], th.
10.12), the theorem is proved.\cqfd

 \th Corollary
\enonce Let $(A,\ell)$ be a general  polarized abelian surface of type
$(1,n)$. The moduli space  $\cM_A(2, \ell,0)$ is smooth irreducible.
\endth

 \th Corollary ([GS])
\enonce Let $ A $ be an abelian surface. The Euler characteristic of
$K_{n-1}(A)$ is $ n^3\sigma(n)$.
\endth

\medskip \centerline{\bf REFERENCES}
\bigskip [AK] A. Altman, S. Kleiman: {\sl Compactifying the Picard
Scheme}, Adv. Math. {\bf 35} (1980), 50--112.

[B1] A. Beauville: {\sl Counting rational curves on K3 surfaces}, e-print
alg-geom 9701019.

[B2] A. Beauville: {\sl Vari\'et\'es k\"ahleriennes dont la premi\`ere
classe de Chern est nulle}, J. Diff. Geom. {\bf 18} (1983), 755--782.

[BL] J. Bryan, N.C. Leung: {\sl The Enumerative Geometry of $K3$ Surfaces
and Modular Forms}, e-print alg-geom 9711031. 

[FL] R. Fahlaoui, Y. Laszlo: {\sl Transform\'ee de Fourier et stabilit\'e
sur les surfaces ab\'eliennes}, Comp. Math. {\bf 79} (1991), 271--278.

[G]  L. G\"ottsche: {\sl A conjectural generating function for numbers of
curves on surfaces},  e-print alg-geom 9711012.

[GS]  L. G\"ottsche, W. Soergel: {\sl Perverse sheaves and the cohomology
of Hilbert schemes of smooth algebraic surfaces},  Math. Ann. 296 (1993),
235--245.

[GH] L. G\"ottsche, D. Huybrechts: {\sl Hodge Numbers of Moduli Spaces of
Stable Bundles on K3 surfaces}, Intern. J. Math. {\bf 7} (1996), 359--372.

[H] D. Huybrechts: {\sl Compact Hyperk\"ahler Manifolds},
Habilitationsschrift, Universit\"at-GH Essen,  1997.

[HL] D. Huybrechts, M. Lehn: {\sl The Geometry of Moduli Spaces of
Sheaves}, Aspects of Mathematics, Vieweg, Braunschweig/Wiesbaden, 1997.

[M1] S. Mukai: {\sl Symplectic structure of the moduli space of sheaves on
an abelian or K3 surface}, Invent. Math.  {\bf 77} (1984), 101--116.

[M2] S. Mukai: {\sl Duality between $D(X)$ and $D(\hat X)$ with its
application to Picard sheaves}, Nagoya Math. J. {\bf 81} (1981), 153--175.

[N] A. Norton: {\sl Analytic moduli of complex vector bundles}, Indiana
Univ. Math. J. {\bf 28} (1979), 365--387.
 
[YZ] S.-T.  Yau, E. Zaslow: {\sl BPS states, string duality, and nodal
curves on K3}, Nuclear Physics B {\bf 471}  (1996), 503--512; also
preprint hep-th 9512121.

\vskip 6mm

\hbox to 72mm{\hrulefill}\parskip=0cm\baselineskip 11pt Olivier {\pc
DEBARRE}

 IRMA -- Math\'ematique --  CNRS 
 
Universit\'e Louis Pasteur

7, rue Ren\'e Descartes

67084 Strasbourg C\'edex -- France

Adresse \'electronique: debarre@math.u-strasbg.fr

\bye